\documentclass{article}
\author{Luis Fernando Garc\'ia-Mora and Hugo Alberto Rincón-Mejía$^1$\
\textsuperscript{1}Corresponding author..}

\title{First elements associated with partial order actions in $\rmod$.}
\date{}

\usepackage{amsmath}
\usepackage{footnote}
\usepackage{todonotes}
\usepackage{multicol}
\usepackage{booktabs}
\usepackage[flushleft]{threeparttable}
\usepackage[font=small,labelfont=bf]{caption}

\usepackage[utf8]{inputenc}
\usepackage{multicol}

\usepackage{blindtext}
\usepackage{amssymb}
\usepackage{amsthm}
\usepackage{amsmath}
\usepackage{dsfont}
\usepackage{mathrsfs}
\usepackage{tikz-cd}
\usetikzlibrary{cd}
\usepackage{comment}

\newcommand{\rpr}{\textbf{$R$-pr}}
\newcommand{\zpr}{\textbf{$\mathbb{Z}$-pr}}
\newcommand{\spr}{\textbf{$S$-pr}}
\newcommand{\trad}{\textbf{$R$-trad}}
\newcommand{\rmod}{R\textit{-Mod}}
\newcommand{\zmod}{\mathbb{Z}\textit{-Mod}}
\newcommand{\rsimp}{R\textit{-simp}}
\newcommand{\rid}{R\textit{-pid}}
\newcommand{\rrad}{R\textit{-rad}}
\newcommand{\rlep}{R\textit{-lep}}
\newcommand{\rler}{R\textit{-ler}}

\newcommand{\submfide}{_{f.i.}\leq_R}

\newcommand{\Mmid}{\hspace{.7mm}\vert\hspace{.7mm}}

\theoremstyle{plain}
\newtheorem{theorem}{Theorem}[section]
\newtheorem{corollary}[theorem]{Corollary}
\newtheorem{lemma}[theorem]{Lemma}

\newtheorem{proposition}[theorem]{Proposition}

\theoremstyle{definition}

\newtheorem{definition}[theorem]{Definition}
\newtheorem{example}[theorem]{Example}

\newtheorem{remark}[theorem]{Remark}

\newtheorem*{explanation}{Explanation}

\begin{document}
\begin{center}

{\Large \bf  First elements of  partial order actions in $\rmod$} 

\vspace{6mm}

{\large \bf Luis Fernando Garc\'ia-Mora and Hugo Alberto Rinc\'on-Mej\'ia.}

\vspace{3mm}
\end{center}

\begin{multicols}{2}

\begin{flushright}\noindent\textbf{MSC2020:}\hspace{1cm}\phantom{.}\\
16S90,16S99,16K99 \end{flushright}

\end{multicols}

\begin{abstract}
	 We explore  concepts of module theory derived  from the notion of primeness, such as first modules, and extend them to more general environments. We also provide descriptions of simple left semiartinian rings, left local rings, semisimple rings, and simple rings in terms of their $\mathscr A$-first modules with respect to a class of preradicals.
\end{abstract}

\section{Introduction}

In this paper, we investigate some concepts of module theory that originated from the idea of primeness, such as first modules,  and we extend them to more general environments. We begin by recalling concepts from module theory as an introduction to the ideas  presented more broadly.\\

In \cite{Johnson}, Johnson calls a module $_RM$  prime if $Ann_R (M)$ = $Ann_R (N)$ for each nonzero submodule $N$ of $M$.  $Ann_R (M)$ denotes the left annihilator of $M$ in $R$, see also \cite{Smi} and \cite{Bea}.\\

When the ring $R$ is commutative, an ideal $I$ is prime if $ab\in I\Rightarrow( (a\in I) \vee (b\in I))$. 
When the ideal $_R0$  is prime,  the ring $R$ becomes an integral domain.
If $R$ is a non commutative ring, an ideal $I$ is prime if $JK\subseteq I\Rightarrow J\subseteq I$  or $ K\subseteq I$.
   So $0$ is prime when $IJ=0 \Rightarrow ( (I=0)   \vee ( J=0)).$ This idea is extended to not necessarily commutative modules and rings by stating that a module $M$ is prime if for any $K$  submodule of $M$ and for any ideal $I$ of $R$, if $IK=0$, then $IM=0$; which is equivalent to Johnson's definition of prime module in \cite{Johnson}.\\

  A class of modules is called pretorsion-free  if it is closed under  taking submodules and forming products. A pretorsion-free class $\zeta$ is considered a torsion-free class if it is closed under extensions. Furthermore, a torsion-free class $\zeta$ becomes a hereditary torsion-free class if it is closed under injective hulls. The collection of hereditary torsion-free classes constitutes a lattice, which we denote as $\mathscr{L}_{\leq, \prod, E, ext}$ for this paper. Similarly, the collection of torsion-free classes forms a big lattice, which we denote as $\mathscr{L}_{\leq, \prod}$.\\
Given a family $\{ \mathcal{A}_i \}_{i \in I}$ of classes in $\mathscr{L}_{\leq, \prod}$, there exists a least member of $\mathscr{L}_{\leq, \prod}$, denoted as $\xi_{\leq, \prod}(\{ \mathcal{A}_i \}_{i \in I})$, such that $\mathcal{A}_j \subseteq \mathscr{L}_{\leq, \prod}$ for each $j \in I$.\\ We write $\xi_{\leq, \prod}(M)$ instead of $\xi_{\leq, \prod}(\{M\})$. \\ Notice that $\xi_{\leq, \prod}(M)$ is the class of left $R$-modules cogenerated by $_RM$.\\

The concept of \textit{prime module} was introduced by L. Bican, P. Jambor, T. Kepka, and P. N\v{e}mec in Section 2 of \cite{Bic1}. A module $M$ is defined as \textit{ prime}  if every nonzero submodule $_RN$ of $_R M$ cogenerates the class  $\xi_{\leq,\prod}(_RM)$. 
\begin{definition}

A module $M$ is  $BJKN$-prime if each nonzero submodule of $M$ cogenerates $_RM$.
 \end{definition}
 
In the same section of \cite{Bic1}, the authors introduce the products of two submodules $A$ and $B$ within a given module $M$.
\begin{definition} Let $A, B$ be two submodules of a module $_RM$. $A \bullet_M B=\sum\{f(A)\vert f\in Hom_R(A,B)\}.$ 
    
\end{definition}

\begin{proposition}[\cite{Bic1} Proposition 2.1]\label{bicprimo}
  Let $M$ be a nonzero $R$-module. $M$ is $BJKN$-prime if and only if, for all proper submodules $L$ and $N$ of $_R M$, the product $L \bullet_M N$ is a nonzero submodule of $M$.
\end{proposition}

   A ring $R$ is a $V$-ring if every simple $R$ module is injective. $R$ is a left semiartinian ring if every nonzero left $R$-module has a nonzero socle. $R$ is a left local ring if two simple left $R$ modules are isomorphic. An $_RM$ module is semisimple if it is the sum of its simple left submodules. $R$ is a semisimple (Artinian) ring if it is semisimple as a $R$-module.\\\vspace{.2cm}

   If $_RM$ is a left module, it is well known that the set of submodules of $M$, which we denote \( \mathscr{L}(_RM),\) is a complete lattice; see Proposition 2.5 in \cite{An}.\\    We say that a submodule $_RN$ of $ _RM$ is fully invariant in $M$,  if $f(N)\leq N$ for all $f\in End(_RM)$. It is easy to see that the set of fully invariant submodules of $ M$ is a complete sublattice of  \( \mathscr{L}(_RM)\). We will denote $_{f.i.}\mathscr{L}(_RM)$, the lattice of fully invariant submodules of $M$, and $N _{f.i.}\leq_R M$ will mean that $N$ is a fully invariant submodule of $M$.

\subsection{Preradicals in $\rmod$.}

We introduce the basic definitions and results of the preradicals in $\rmod$. For more information on preradicals, see \cite{Bic}, \cite{Clark}, \cite{pre}, \cite{preII}, \cite{preIII}, and \cite{Es}.
A preradical $\sigma$ on $\rmod$ is a functor $\sigma:\rmod\rightarrow \rmod$ such that:
\begin{enumerate}
    \item $\sigma(M)\leq M$ for each $M\in \rmod$.
    \item For each $R$-morphism $f:M\rightarrow N$,  the following diagram is commutative: \begin{center}
\begin{tikzcd}
 M \ar{r}{f}& N\\
\sigma(M)\ar[hookrightarrow]{u}\ar{r}{f_\downharpoonright^\upharpoonright} & \sigma(N).\ar[hookrightarrow]{u}
\end{tikzcd}
\end{center}
\end{enumerate}
Recall that for each $\beta\in \rpr$ and each $\{M_i\}_{i\in I}$ family of $R$-modules, we have $\beta (\bigoplus \limits _{i\in I}M_i)=\bigoplus \limits _{i\in I}\beta (M_i)$, see Proposition I.1.2 of \cite{Bic}. 
$\rpr$ denotes the collection of all preradicals in $\rmod$. By Theorem 7 of \cite{pre}, $\rpr$ is a big lattice with smallest element $\underline{0}$ and largest element $\underline{1}$, where:\begin{enumerate}
    \item the order in  $\rpr$ is given by $\alpha\preceq \beta$ if $\alpha(M)\leq \beta(M)$ for every $M\in \rmod$;
    \item for any family of preradicals $\{\sigma_i\}_{i\in I}$ in $\rmod$ the supremum and infimum of the family are given, respectively, by \begin{enumerate}
        \item $(\bigvee\limits_{i\in I}\sigma_i)(M)=\sum\limits_{i\in I}\sigma_i(M)$ and 
        \item $(\bigwedge\limits_{i\in I}\sigma_i)(M)=\bigcap\limits_{i\in I}\sigma_i(M)$.
\end{enumerate}
\end{enumerate}  
Recall that $\sigma\in \rpr$ is idempotent if $\sigma \circ \sigma =\sigma$.  $\sigma$ is radical if $\sigma(M/\sigma(M))=0$, for each $M\in \rmod$.  $\sigma$ is a left exact preradical if it is a left exact functor. 
$\sigma$ is $t$-radical (also known as a cohereditary preradical) if $\sigma(M)=M\sigma(R)$. 
Recall that $\sigma$ is a $t$-radical if and only if $\sigma$ preserves epimorphisms. See Exercise 5  of Chapter VI \cite{Es},  and that  $\sigma$ is a left exact preradical if and only if, for each submodule $N$ of a module $M$, we have $\sigma(N)=\sigma(M)\cap N$; see Proposition 1.7 of Chapter VI of \cite{Es}. 
We will denote $\rid$, $\rrad$, $\rlep$, $\trad$ and $\rler$, the collections of idempotent preradicals, radicals, left exact preradicals, $t$-radicals, and left exact radicals, respectively. For a $\sigma \in \rpr$ we will denote by $\widehat{\sigma}$ the largest idempotent preradical smaller than $\sigma$ and by $\overline{\sigma}$ the least radical greater than $\sigma$. See Proposition 1.5 in Chapter VI of \cite{Es}.

\subsection{Relationships between preradicals and fully invariant submodules.}

Recall for a $_RN$ fully invariant submodule of $_RM$, the preradicals $\alpha_N^M$ and $\omega_N^M$ defined in \cite{pre}, Definition 4.
\begin{definition}
Let $_RN$ be a fully invariant submodule of $_RM$ and $U\in \rmod$. 
\begin{enumerate}
\item $\alpha_N^M(U)=\sum\lbrace f(N) \Mmid f\in Hom_R(M,U) \rbrace$ and 
\item $\omega_N^M(U)=\bigcap\lbrace f^{-1}(N) \Mmid f\in Hom_R(U,M) \rbrace$.
\end{enumerate}  
\end{definition}

More generally, for $ _RN$ a submodule of $_RM$, we can define $\beta_N^M$  the preradical such that $\beta_N^M(U)=\sum\lbrace f(N) \Mmid f\in Hom_R(M,U) \rbrace$. Note that $A\bullet_MB=\beta_A^M(B)$ for each $A,B\leq M$.

\begin{remark}
  Note that for each $M\in\rmod$, we have that $(\omega_0^M:\omega_0^M)(M)/\omega_0^M(M)=\omega_0^M(M/\omega_0^M(M))=\omega_0^M(M/0)=\omega_0^M(M)=0$. Therefore, $\omega_0^M$ is a radical.
\end{remark}
\begin{remark}
Let $M\in \rmod$ and let $N$ be a fully invariant submodule of $M$. By Proposition 5 of \cite{pre},  $\alpha_N^M :\rmod\rightarrow \rmod$ is the least preradical $\rho$ such that $\rho(M)=N$, and $\omega_N^M :\rmod\rightarrow \rmod$ is the largest preradical $\rho$  such that $\rho(M)=N$. It is easy to see that $\{\rho\in \rpr\mid \rho(M)=N\}=[\alpha_N^M, \omega_N^M]$, an interval in $R$-pr. 
\end{remark}

We can see that $\alpha_ M^M(L)$ is the trace of $M$ in $L$, that is, $\alpha_M^M=tr_{M}$.
 Furthermore, $\omega^M_0(L)$ is the
reject of $M$ in $L,$ and is the smallest submodule of $L$ such that $L/\omega^M_0(L)$
embeds in a product of copies of $M$.

\begin{remark}\label{nr1}
    If  $L\submfide K\submfide M$,   then $L\submfide M$.
\end{remark}

\begin{proof}
If $f:M\rightarrow M$  is an $R$-morphism,  then $f$ restricts to an $R$-endomorphism
of $K$, which in turn restricts to an endomorphism of $L$, see \cite{Clark} Exercises 6.31.(1) and (2).
\end{proof}

As each simple $R$-module is isomorphic to an $R/I$, where $I$ is a maximal left ideal of $R$, a set of representatives of the isomorphism classes of simple modules $ \rsimp$ can be chosen.
Recall that $soc:\rmod\rightarrow \rmod$ is an idempotent preradical such that $soc(M)$ is the largest semisimple submodule of $M$. Furthermore, $soc=\vee \{\alpha_S^S\hspace{1mm}\mid \hspace{1mm} S\in R\textit{-simp }\}$.

\section{$\mathscr{P}$-first elements in a lattice.}

\begin{definition}
A lattice $\mathscr{L}$ is bounded if elements $0,1\in \mathscr{L}$ exist such that $0\leq x\leq 1$ for each $x\in \mathscr{L}$. 
\end{definition}

The following definition appears in \cite{Abu}.
\begin{definition}
Let $\mathscr{L}=(L,\leq ,\vee, \wedge)$ be a lattice and let $\mathscr{P}=(P,\leq ' )$ be a poset. A $\mathscr{P}$-action on $\mathscr{L}$ is a function $\rightharpoonup: \mathscr{P}\times \mathscr{L}\rightarrow \mathscr{L}$ satisfying the following conditions for all  $s,t\in  P$ and $x,y\in L$:
\begin{itemize}
\item[i)] $s\leq ' t$ $\Rightarrow$ $s\rightharpoonup x \leq t\rightharpoonup x$.
\item[ii)] $x\leq  y$ $\Rightarrow$ $s\rightharpoonup x \leq s\rightharpoonup y$.
\item[iii)] $s\rightharpoonup x\leq x$.
\end{itemize}
\end{definition}

Since lattices are specific types of posets,  we can explore their actions on other  lattices as specific cases of the above definition. 
We consider especially the lattices $\mathscr{L}(_RM)$, $\rpr$ and   $_{fi}\mathscr{L}(_RM)$ acting  on other partial order sets or other partial sets acting on them.
We distinguish $\mathscr{L}(_RR)$, the lattice of left ideals of $R$,  from  $Lat(R)$,  which will denote the lattice of two-sided ideals of $\rmod$.
\begin{example}\label{1.1}
The lattice $\mathscr{L}(_RR)$ acts on $\mathscr{L}(_RM)$. 
 \begin{center}. '
           \begin{tikzcd}[row sep=tiny]
               \mathscr{L}(_RR)\times \mathscr{L}(_RM)\arrow[r,"\rightharpoonup"]& \mathscr{L}(_RM) \\
                (I,N)\arrow[r, mapsto] & IN.
           \end{tikzcd}
        \end{center}
        
\end{example}

\begin{example}\label{1}
The lattice $\rpr$ acts on $\mathscr{L}(_RM)$. 
 \begin{center}
           \begin{tikzcd}[row sep=tiny]
               \rpr\times \mathscr{L}(_RM)\arrow[r,"\rightharpoonup"]& \mathscr{L}(_RM) \\
                (r,N)\arrow[r, mapsto] & r(N).
           \end{tikzcd}
        \end{center}
        
\end{example}

\begin{remark}
Let $\mathscr{P}=(P,\leq ' )$ be a poset, $\mathscr{L}=(L, \leq, \wedge, \vee,0, 1)$ be  a bounded lattice with a  $\mathscr{P}$-action  
$\rightharpoonup:\mathscr{P}\times \mathscr{L}\rightarrow \mathscr{L}$ and $\mathscr{Q}$ be a subposet of $\mathscr{P}$. Then, the restriction of $\rightharpoonup$ to $\mathscr{Q}\times \mathscr{L}$ is a $\mathscr{Q}$-action.
\end{remark}

\begin{definition}\label{2}
Let $\mathscr{P}=(P,\leq ' )$ be a poset, $\mathscr{L}=(L, \leq, \wedge, \vee,0, 1)$ be  a bounded lattice with a $\mathscr{P}$-action  
$\rightharpoonup:\mathscr{P}\times \mathscr{L}\rightarrow \mathscr{L}$ and $x\in L\setminus \{0\}$. We say that:
\begin{itemize}
\item[i)] $x$ is $\mathscr{P}$-first if  for each $0 \neq  z\leq x$ and  $s\in P $\begin{center}
    $s\rightharpoonup z=0$ $\Rightarrow$ $s\rightharpoonup x=0$.
\end{center}   

\item[ii)] $x$ is  $\mathscr{P}$-prime if   for each $  z\in L$ and  $s\in P $
\begin{center}
    $s\rightharpoonup z\leq x$ $\Rightarrow$ ($s\rightharpoonup 1\leq x$ or  $z\leq x$).
\end{center}

\end{itemize}

\end{definition}

\begin{remark}
 When we apply these definitions to the action of   Example \ref{1.1},  we note that the $\mathscr{L}(R)$-first modules are precisely the prime modules. 
\end{remark}

We omit the easy proof of the following proposition.

\begin{proposition}\label{f0}
Each order-preserving function $f:\mathscr{P}\rightarrow \mathscr{Q}$  induces  a correspondence between the class of  $\mathscr{Q}$-actions on $\mathscr{L}$ and the class of $\mathscr{P}$-actions on $\mathscr{L}$. Explicitly,  \begin{tikzcd}[row sep=tiny]
\mathscr{Q}\times \mathscr{L}\overset{\rightharpoonup}{ \longrightarrow }\mathscr{L}\arrow[r,mapsto, "(\_)_f"]& \mathscr{P}\times \mathscr{L}\overset{\rightharpoonup_f}{ \longrightarrow}\mathscr{L}
           \end{tikzcd}, where for each $a\in \mathscr{P}$ and $x\in \mathscr{L}$,  $a\rightharpoonup_fx:=f(a)\rightharpoonup x$.  
 
\end{proposition}

\begin{lemma}\label{f.5}
Let $\mathscr{P}=(P,\leq ' )$ be a poset, $\mathscr{L}=(L, \leq, \wedge, \vee,0, 1)$ be  a bounded lattice with a $\mathscr{P}$-action  
$\rightharpoonup:\mathscr{P}\times \mathscr{L}\rightarrow \mathscr{L}$ and let $f:\mathscr{P}\rightarrow \mathscr{Q}$ be  an  order preserving function. Then,  $\mathscr{Q}$-first elements of $\mathscr{L}$  respect to $\rightharpoonup$,  are $\mathscr{P}$-first elements of $\mathscr{L}$ respect to $\rightharpoonup_f$.
\end{lemma}

In particular, note that we can restrict an $\mathscr{Q}$-action in $\mathscr{L}$ to a subset $\mathscr{P}$ of $\mathscr{Q}$.\\

It is easy to verify the following corollary.

\begin{corollary}\label{f3}
Let $\mathscr{P}=(P,\leq ' )$ be a poset, $\mathscr{L}=(L, \leq, \wedge, \vee,0, 1)$ be  a bounded lattice with a $\mathscr{P}$-action  
$\rightharpoonup:\mathscr{P}\times \mathscr{L}\rightarrow \mathscr{L}$ and let $\mathscr{Q}$ be a subset  of $\mathscr{P}$. Then any  $\mathscr{P}$-first element of $\mathscr{L}$ is $\mathscr{Q}$-first.
\end{corollary}

\begin{definition}
    Let $\mathscr{L}=(L, \leq, \wedge, \vee,0, 1)$ be  a bounded lattice and $y,x\in L$ such that $y\leq x$. We denote $[y,x]$ the bounded lattice of all $z\in L$ such that $y\leq z \leq x$.
\end{definition}

\begin{remark}
    Let $\mathscr{P}=(P,\leq ' )$ be a poset, $\mathscr{L}=(L, \leq, \wedge, \vee,0, 1)$ be  a bounded lattice with a $\mathscr{P}$-action  
$\rightharpoonup:\mathscr{P}\times \mathscr{L}\rightarrow \mathscr{L}$ and $x\in L$. Then $\rightharpoonup_\downharpoonright:\mathscr{P}\times [0,x]\rightarrow [0,x]$ is well defined and is a $\mathscr{P}$-action.
\end{remark}
Recall Remarks 1.8(6) and 1.8(7) of \cite{Abu}.
\begin{remark}\label{primfrt}
    Let $\mathscr{P}=(P,\leq ' )$ be a poset, $\mathscr{L}=(L, \leq, \wedge, \vee,0, 1)$ be  a bounded lattice with a $\mathscr{P}$-action  
$\rightharpoonup:\mathscr{P}\times \mathscr{L}\rightarrow \mathscr{L}$ and $x\in L\setminus\{0\}$. Then $x$ is $\mathscr{P}$-first in $\mathscr{L}$ if and only if $0$ is $\mathscr{P}$-prime in $[0,x]$.
\end{remark}

 Remark \ref{primfrt} gives us a way to relate the concept of $\mathscr{P}$-first  with the concept of $\mathscr{P}$-prime.

\begin{lemma}\label{f2}
Let $\mathscr{P}=(P,\leq ' )$ be a poset, $\mathscr{L}=(L, \leq, \wedge, \vee,0, 1)$ be  a bounded lattice with a  $\mathscr{P}$-action  
$\rightharpoonup:\mathscr{P}\times \mathscr{L}\rightarrow \mathscr{L}$. All atoms of $\mathscr{L}$ are $\mathscr{P}$-first.
\end{lemma}

\begin{proof}
Let $a$ be an atom of $\mathscr{L}$. Let $s\in \mathscr{P}$ and $ z\leq x$.  If $s\rightharpoonup z=0$ and $z\neq 0$, then $a=z$. From this it follows that $s\rightharpoonup a=0$. Therefore, $a$ is $\mathscr{P}$-first.
\end{proof}

\section{\rpr-first modules.}
We can apply Definition \ref{2}  to the action described in   Example \ref{1} to define first modules with respect to preradicals.

\begin{definition}\label{defsec}
Let $M\in \rmod$ and $N\in \mathscr{L}(M)\setminus \{0\}$. We say that $N$ is  a  \textbf{$R$-pr}-first submodule of $M$ if $\alpha(K)=0$ and $K\leq N$ imply $\alpha(N)=0$,  for all $K\in \mathscr{L}(M)\setminus \{0\}$ and all $\alpha\in \rpr $.
\end{definition}

\begin{example}
The nonzero submodules of $_\mathbb{Z}\mathbb{Z}$ are $\zpr$-first submodules of $\mathbb{Z}$, since for each $N,K\in\mathscr{L}(\mathbb{Z})\setminus \{0\}$ we have  $N\cong K$.
\end{example}

\begin{lemma}
Let $M\in \rmod$. Every simple submodule of $M$ is $\rpr$-first.
\end{lemma}

\begin{proof}
This follows from Lemma \ref{f2}.
\end{proof}

The following lemma is a consequence of  Lemma 9.2 of \cite{An}, and of the fact that for any $S,M\in \rmod$ with $S$ simple, $\alpha_S^S(M)$ is a semisimple homogeneous $R$-module.

\begin{lemma}\label{5}
Let $M\in \rmod$ and $N$ be a semisimple nonzero submodule of $M$. Then  $N$ is $\rpr$-first submodule of $M$  iff  $N\cong S^{(I)}$ for some set $I$ and some simple $R$-module $S$.

\begin{proof}
 $\Rightarrow)$ If $0 \neq L \leq N \leq M$ and $\sigma(L) \neq 0$, then since $L = soc(L) = \oplus \{S_i \vert i \in I\}$ for some set $I$, and a preradical commutes with coproducts, we have $0 \neq \sigma(L) = \oplus \{ \sigma(S_i) \vert i \in I\}$. Consequently, $S_i = \sigma(S_i) \neq 0$ for each $i \in I$. Thus, $0 \neq \oplus \{ \sigma(S_i) \vert i \in I\} \leq \sigma(N).$\\$\Leftarrow)$  If  $N\cong S^{(I)}$ and $\sigma(N)\neq0$, as above,  $\sigma(S)=S$. If $0\neq L\leq N$, then $L=S^{(J)} $, for some set $J$, then $\sigma(L)=L\neq 0.$
\end{proof}
\end{lemma}

\begin{remark}
We find that if $N$ is a submodule of $M$, then $N$ is $\rpr$-first  in $\mathscr{L}(M)$ if and only if $N$ is $\rpr$-first in $\mathscr{L}(N)$. In this sense, we will say that an $R$-module $M$ is $\rpr$-first  if it is $\rpr$-first as a submodule of itself.
\end{remark}

\begin{example}
$_\mathbb{Z}\mathbb{Z}$ is $\zpr$-first, since for any $N\in\mathscr{L}(\mathbb{Z})\setminus \{0\}$, we have $N\cong \mathbb{Z}$.
\end{example}

Recall Definition 42 of \cite{primepre}.
\begin{definition}
An $R$-module $M$ is diuniform if any fully invariant submodule $N$ of $M$ is essential in $M$.
\end{definition}

\begin{proposition}\label{6}
Let $M$ be an $\rpr$-first $R$-module. Then $M$ is diuniform.
\end{proposition}

\begin{proof}
Let $_RM$ be $\rpr$- first, $N\in  \mathscr{L}_{fi}(M)\setminus \{ 0\}$ and $K\in  \mathscr{L}(M)\setminus \{ 0\}$. We have $\alpha_N^M(M)=N\neq 0$ and $M$  is $\rpr$-first, so $\alpha_N^M(K)\neq 0$.  Moreover, $\alpha_N^M(K)\leq N$ and $\alpha_N^M(K)\leq K$, since $\alpha_N^M$ is a preradical, and $N$ is fully invariant, so $0\neq \alpha_N^M(K)\leq N\cap K$, then $ N\cap K\neq 0$. Therefore, $N$ is essential in $ \mathscr{L}(M)$.
\end{proof}

The proposition \ref{6} shows that $\mathscr{L}_{f.i.}(M)$ is a uniform lattice if $M$ is $\rpr$-first. The following example shows that the converse of Proposition \ref{6} is not true.

\begin{example}\label{7}
We have $Hom_\mathbb{Z}(\mathbb{Z}_{p^\infty},\mathbb{Z}_p)=0$, so $\alpha_{\mathbb{Z}_p}^{\mathbb{Z}_{p^\infty}}(\mathbb{Z}_p)=0$ and $\alpha_{\mathbb{Z}_p}^{\mathbb{Z}_{p^\infty}}(\mathbb{Z}_{p^\infty})=\mathbb{Z}_p\neq 0$, so $_\mathbb{Z}\mathbb{Z}_{p^\infty}$ is not $\rpr$-first.
\end{example}

The following proposition characterizes the $\rpr$-first modules.  

\begin{proposition}\label{pf}
The following statements are equivalent:

\begin{itemize}
\item[(1)] $_RM$ is $\rpr$-first.
\item[(2)] $M$ is $BJKN$-prime.
\item[(3)] Each nonzero  cyclic submodule of $M$ cogenerates $M$.
\item[(4)] For any $x,y\in M\setminus\left\{ 0\right\} $ $\exists f_x:M\longrightarrow Ry$
such that $f_x\left(x\right)\neq 0$.

\end{itemize}
\end{proposition}
\begin{proof}
\begin{itemize}
\item[] 
    \item[$(1)\Rightarrow (2)$]  Suppose that $A,B$ are nonzero submodules of $M.$ Let $\beta$  preradical defined by $\beta\left(N\right)=\sum\left\{ f\left(A\right)\mid f:M\longrightarrow N\right\} .$
    Then $\beta\left(M\right)\geqq A\neq 0,$ so by hypothesis, $\beta\left(B\right)\neq 0.$
Thus $\exists f:M\longrightarrow B$ such that $f\left(A\right)\neq0.$
This means that $A\bullet B\neq 0$, by Proposition \ref{bicprimo}.
\item[$(2)\Rightarrow (3)$] It is inmediate.

    \item[$(3)\Rightarrow (4)$]  Let $x,y\in M\diagdown\left\{ 0\right\} $, as $Ry$ cogenerates $M,$  there exists a monomorphism $\varphi:M\longrightarrow\left(Ry\right)^{Z},$ for some set $Z.$ Then $\exists z\in Z$ such that for the  canonical projection  $\pi_{z},$ $\left(\pi_{z}\circ\varphi\right)\left(x\right)\neq 0.$ 
    \item[$(4)\Rightarrow (1)$]  Suppose $\sigma$ is a preradical such that $\sigma\left(M\right)\neq0.$ Let $B\neq0$ be a submodule of $M.$
If  $0\neq x\in\sigma\left(M\right)$ and  $0\neq y\in B,$  there exists $f:M\longrightarrow Ry$
such that  $f\left(x\right)\neq0.$ Then, the following diagram is commutative:
\begin{center}
\begin{tikzcd}
 M \ar{r}{f}& Ry \\
\sigma(M)\ar[hookrightarrow]{u}\ar{r}{f\vert} & \sigma(Ry)\ar[hookrightarrow]{u}\\
x\ar[mapsto]{r} & f(x).
\end{tikzcd} 

\end{center}
Thus, $0\neq\sigma\left(Ry\right)\leq\sigma\left(B\right)$.
\end{itemize}
\end{proof}

\section{First modules relative to a subclass of $\rpr$.}
Let $\mathscr{A}\subseteq \rpr$. We can consider the action of $\mathscr{A}$ on the submodule lattice of a module induced by the inclusion of $\mathscr{A}$ in $\rpr$ and the $\rpr$ action described in Example \ref{1}, as in Lemma \ref{f0}, for the following definition.

\begin{definition}
Let $ M\in \rmod$ and $\mathscr{A}\subseteq \rpr$. We say that:
\begin{itemize}
\item[i)] $M$ is  $\mathscr{A}$-first if $0\neq M$ and $\alpha(K)=0$ imply $\alpha(M)=0$,  for all $K\in \mathscr{L}(M)\setminus \{0\}$, and for all $\alpha\in \mathscr{A} $.
\item[ii)] $M$ is  $\mathscr{A}$-fully first if $\alpha(K)\neq 0$,  for all $K\in \mathscr{L}(M)\setminus \{0\}$ and for all $ \alpha\in \mathscr{A} $.

\end{itemize}
\end{definition}

Notice that each $t$-radical $\sigma$  can be written   as $ \sigma (R)\cdot\_$ and recall that for each  preradical  $\sigma$, $\sigma(R)$ 
is a two-sided ideal of $R$. Thus, there is a lattice isomorphism between
$\trad$ and Lat(R) sending $\sigma$ to $\sigma(R)$. The following proposition follows from Lemma \ref{f.5}.

\begin{proposition}
For each $M\in \rmod$ the following statements are equivalent.\begin{itemize}
    \item[(1)] $M$ is a prime module.
    \item[(2)] $M$ is a $\trad$-first module.
\end{itemize} 
\end{proposition}

\begin{theorem}\label{15}
For a ring $R$,  the following statements are equivalent:
\begin{itemize}
    \item[(1)] $R$ is a simple ring.
    
     \item[(2)] Any $R$-module is prime.
\end{itemize}
\end{theorem}

\begin{proof}
\begin{itemize}
    \item[]
    
   \item[(2)$\Rightarrow$ (1)] It is clear.
   \item[(1)$\Rightarrow$ (2)] Let $I$ be a proper ideal of $R$. By hypothesis, $R/I\oplus R$ is prime,  and $I(R/I)=0$, so $I(R/I)\oplus IR=I(R/I\oplus R)=0$, then $I=IR=0$.\\
   Therefore, $R$ is simple.
\end{itemize}
\end{proof}

\begin{example}\label{ejem}
Let $R$ be the ring of linear endomorphisms  of the  $\mathbb{R}$ vector space $\mathbb{R}^{(\mathbb{N})}$. 
Consider $I=\{f \in R\Mmid rank(f)\in \mathbb{N}\}$ and $S=R/I$. As $S$ is a simple ring, any $S$-module is first by Theorem \ref{15}.\\
On the other hand, $S$ is not a left semisimple ring, so $ \mathbb{P}_{\spr}\neq \rmod$ by Theorem \ref{14.3}.
As $ \mathbb{P}_{\spr}\neq \rmod$, there exists $M\in \rmod$ such that $M$ is prime, but $M$ is not $\spr$ first. Therefore, a module $_RM$ can be prime without being $\rpr$ first ($BJKN$-prime).
\end{example}

\begin{proposition}\label{rpid1}
    Let $M\in \rmod$ be a nonzero module. The following statements are equivalent.
    \begin{itemize}
        \item[(1)] $M$ is a $\rid$-first module.
        \item[(2)] For each nonzero submodules  $N$ and  $K$ of $M$, there exists a nonzero morphism $f:N\rightarrow K$.
    \end{itemize}
\end{proposition}

\begin{proof}
    \begin{itemize}
        \item[]
        \item[$(1)\Rightarrow (2) $] Let $N, K\in \mathscr{L}(M)\setminus \{0\}$. Then $tr_N\in \rid$ and $0\neq N\leq tr_N(M)$. Moreover, $M$ is $\rid$-first and $0\neq K\leq M$, so $tr_N(K)\neq 0$, which implies that  there exists a nonzero morphism $f:N\rightarrow K$.
        \item[$(2)\Rightarrow (1)$] Let $N\in \mathscr{L}(M)\setminus \{0\}$ and $\sigma \in \rid$ such that $\sigma (M)\neq 0$. Then there exists a nonzero morphism $f:\sigma(M)\rightarrow N$. As  $\sigma\in \rid$, the following diagram is commutative
        \begin{center}
\begin{tikzcd}
 \sigma(M) \ar{r}{f}& N \\
\sigma(\sigma(M))\ar[hookrightarrow]{u}\ar{r}{f\vert} & \sigma(N).\ar[hookrightarrow]{u}
\end{tikzcd} 

\end{center}
Thus, $0\neq f(\sigma(M))=f(\sigma(\sigma(M)))\leq\sigma(N)$. Therefore, $M$ is $\rid$-first.
\end{itemize}
\end{proof}

\begin{proposition}\label{endoprimo1}
    Let $M\in \rmod$ be a nonzero module. If $M$ is retractable and $End_R(M)$ is prime, then $M$ is a $\rid$-first module.
\end{proposition}

\begin{proof}
    Let $\sigma\in\rid$ such that $\sigma(M)\neq 0$. Let $N$ be a nonzero submodule of $M$. Consider the following.
    \begin{eqnarray*}
        I&=&\{f\in End_R(M)\mid f(M)\subseteq \sigma(M)\},\\
        X&=&\{f\in End_R(M)\mid f(M)\subseteq N\} \textit{ and }\\
        J&=&\{\sum\limits_{i=1}^ns_if_i\mid n\in \mathbb{N}, s_i\in End_R(M) \textit{ and } f_i\in X\}.
    \end{eqnarray*}

Is clear that $I$ and $J$ are two ideals of $End_R(M)$. Since $M$ is retractable, there exist two nonzero morphisms $g : M \to N$ and $f : M \to \sigma(M)$. So, if $\overline{g}:M\to M$ and $\overline{f}:M\to M$ are the morphisms defined by the diagrams  \begin{tikzcd}
    M\ar[r, "g"']\ar[rr, bend left, "\overline{g}"] & N\ar[r, hookrightarrow] & M
    \end{tikzcd} and \begin{tikzcd}
    M\ar[r, "f"']\ar[rr, bend left, "\overline{f}"] & \sigma(M)\ar[r, hookrightarrow] & M
    \end{tikzcd}, then $0\neq \overline{f}\in I $ and $0\neq \overline{g}\in X\subseteq J$, so $I$ and $J$ are two nonzero ideals of $End_R(M)$. Now, since $End_R(M)$ is prime, then  $JI\neq 0$.\\
 There then exist $h\in End_R(M)$, $l\in X$ and $k\in I$ such that $hlk\neq 0$, thus $lk\neq 0 $. Moreover,  $0\neq Im(lk)=l(Im(k))\subseteq l(\sigma(M))$.  Since $l(\sigma(M))\subseteq l(M)\subseteq N$, the restriction $\overline{l}=l\vert_{\sigma(M)}^N:\sigma(M)\to N$ is a nonzero morphism. As  $\sigma\in \rid$, the following diagram is commutative:
    \begin{center}
    \begin{tikzcd}
 \sigma(M) \ar{r}{\overline{l}}& N \\
\sigma(M)=\sigma(\sigma(M))\ar[hookrightarrow]{u}\ar{r}{\overline{l}\vert} & \sigma(N).\ar[hookrightarrow]{u}
\end{tikzcd} 
\end{center}
Thus, $0\neq \overline{l}(\sigma(M))=\overline{l}(\sigma(\sigma(M)))\leq\sigma(N)$. Therefore, $M$ is $\rid$-first.
\end{proof}

\begin{example}
    Let $R=\mathbb{Z}_{p^2}=\mathbb{Z}/(p^2\mathbb{Z})$. $_RR$ has only 3 submodules $R$, $\langle p+p^2\mathbb{Z}\rangle$ and the zero submodule. So it is clear that $R$ is a  $\rid$-first module by the Proposition \ref{rpid1}. On the other hand, $End(_RR) \cong R$, and $R$ is not a prime ring. Therefore, the converse of Proposition \ref{endoprimo1} does not hold.
\end{example}


\begin{proposition}
For a nonzero $_RM$, the following statements are equivalent:

\begin{itemize}
\item[(1)] $M$ is $\rrad$-first.
\item[(2)] $M$ is $BJKN$-prime.
\end{itemize}
\end{proposition}

\begin{proof}
    \begin{itemize}
        \item[]
        \item[$(1)\Rightarrow (2) $] Let $N\in \mathscr{L}(M)\setminus \{0\}$. Then $\omega^N_0\in \rrad$ and $\omega_0^N(N)=0$. Since $M$ is $\rrad$-first and $0\neq N\leq M$, then $\omega_0^N(M)=0$, which implies that  $M$ is cogenerated by $N$. Therefore, $M$ is $BJKN$-prime.
        \item[$(2)\Rightarrow (1)$] Suppose $M$ is a $BJKN$-prime, then, by Proposition \ref{pf}, $M$ is $\rpr$-first, and by Lemma \ref{f.5},  $M$ is $\rrad$-first.
\end{itemize}
\end{proof}

\section{The classes of $\mathscr{A}$-first modules and the class of $\mathscr{A}$-fully first modules}

Let $\sigma \in \rpr$. $\mathbb{T}_\sigma=\{ M\in \rmod \mid \sigma(M)=M\}$ is the class of $\sigma$-pretorsion modules. $\mathbb{F}_\sigma=\{ M\in \rmod \mid \sigma(M)=0\}$ is the class of $\sigma$-pretorsion free  modules. For $\mathscr{A}\subseteq \rpr$, we define  $\mathbb{T}_{\mathscr{A}}:=\bigcap\limits_{r\in \mathscr{A}}\mathbb{T}_{r}$ and 
 $\mathbb{F}_{\mathscr{A}}:=\bigcap\limits_{r\in \mathscr{A}}\mathbb{F}_{r}$.\\ Now, we  denote 
\begin{itemize}
    \item $\mathscr{P}_{\mathscr{A}}$ is the class of all $\mathscr{A}$-fully first modules;
    \item $\mathbb{P}_{\mathscr{A}}$ is the class of all $\mathscr{A}$-first modules and the zero module.
\end{itemize}
We will use  $\mathscr{P}_{\sigma}$ and  $ \mathbb{P}_{\sigma}$ instead of  $\mathscr{P}_{\{\sigma\}}$  and $ \mathbb{P}_{\{\sigma\}},$  respectively.

\begin{remark}\label{8}
For any $\sigma,\beta\in \rpr$ and $\mathscr{A},\mathscr{B}\subseteq \rpr$, it follows that:
\begin{itemize}
\item[(1)]  $\mathscr{P}_{\mathscr{A}}=\bigcap\limits_{r\in \mathscr{A}}\mathscr{P}_{r}$ and 
 $\mathbb{P}_{\mathscr{A}}=\bigcap\limits_{r\in \mathscr{A}}\mathbb{P}_{r}$.
 \item[(2)] $\mathbb{P}_{\sigma}=\mathscr{P}_{\sigma} \cup \mathbb{F}_{\sigma}$.
 \item[(3)]  $\mathscr{P}_{\mathscr{A}}\subseteq \mathbb{P}_{\mathscr{A}}$ and 
 $\mathbb{F}_{\mathscr{A}}\subseteq \mathbb{P}_{\mathscr{A}}$.
 \item[(4)] If $\sigma \preceq \beta$, then $\mathscr{P}_\sigma\subseteq \mathscr{P}_\beta$.
 \item[(5)] If $\mathscr{B}\subseteq \mathscr{A}$ then $\mathscr{P}_\mathscr{A}\subseteq \mathscr{P}_\mathscr{B}$ and $\mathbb{P}_\mathscr{A}\subseteq \mathbb{P}_\mathscr{B}$.
\end{itemize}
\end{remark}

To see the point $(2)$ of Remark \ref{8}, let us notice that $M\in\mathbb{P}_{\sigma}$ if and only if for each $0\neq N\leq M$ we have $\sigma(N)= 0 \implies \sigma (M)= 0$. That is equivalent to $\sigma(M)=0$ or $\sigma(N) \neq 0$ for each $0\neq N\leq M$. 

\begin{remark}
If $N$ is a nonzero element of $\mathscr{L}_{fi}(M)$, then $N \leq_e M$. This is because if $M$ belongs to $\mathscr{P}_{\alpha_N^M}$, then $\alpha_N^M(M) = N \neq 0$. Consequently, for any nonzero $K$ such that $K \leq M$, $\alpha_N^M(K) \neq 0$. Since $\alpha_N^M(K) \subseteq N \cap K$, it follows that $N \leq_e M$.
\end{remark}

\begin{example}\label{ExampleZclass}
\textbf{$\alpha_{\mathbb{Z}_p}^{\mathbb{Z}_p}$-first modules in $\mathbb{Z}$-mod.}\\
Let $p\in \mathbb{Z}$ be a prime number,  and let  $\sigma=\alpha_{\mathbb{Z}_p}^{\mathbb{Z}_p}$ .  For a non zero $\mathbb{Z}$-module $M$ the following hold:
\begin{itemize}
    \item[(1)]  $M$ is $\sigma$-fully first if and only if $M$ is a $p$-group.\\
    If $N\in \mathscr{P}_{\sigma}$ and $0\neq x\in N$, then $\sigma(\langle x \rangle)\neq 0$, so that there exists $0\neq f\in Hom_\mathbb{Z}(\mathbb{Z}_p, N)$. From this we find that, for some $r\in \mathbb{Z}$,  $0\neq f(1+p\mathbb{Z})=x^r$  is an element of order $p$ of $N$. We conclude that $x$ has an order that is a multiple of $p$. If $o(x)=p^ns$ with $n\in \mathbb{N}$ and $p$ does not divide $s$, then $o(x^{p^n})=s$, so $Hom_\mathbb{Z}( \mathbb{Z}_p, \langle x^{p^n} \rangle)=0$. So, $\sigma (\langle x^{p^n} \rangle)=0$, which implies that $\langle x^{p^n} \rangle=0$.  Therefore, $o(z)=p^n$, so it follows that $N$ is a $p$ group.
    Conversely , if $N$ is a $p$-group and $0\neq K\leq N$, then there exists $x\in K\setminus \{0\}$, so there exists $n\in \mathbb{N}\setminus \{0\}$ such that $o(x)=p^n$. Moreover,  $o(x^{p^n-1})=p$, which implies $\langle x^{p^n-1} \rangle\cong \mathbb{Z}_p$ and thus $Hom_\mathbb{Z}( \mathbb{Z}_p, K)\neq 0$, from which it follows that $\sigma (K)\neq 0$. Consequently, $N\in \mathscr{P}_{\sigma}$ if and only if $N$ is a $p$-group.\\
    \item[(2)] $M$ is $\sigma$-first if and only if $M$ has no elements of order $p$ or $M$ is a $p$-group.

If $N\in \zmod$ has an element $a$ of order $p$, then $K=\langle a\rangle \cong \mathbb{Z}_p $, so $\sigma (N)\neq 0$. On the other hand, if $\sigma (N)\neq 0$, then there exists $0\neq f\in Hom_\mathbb{Z}(\mathbb{Z}_p, N)$, so $f(1+p\mathbb{Z})$ is an element of order $p$ of $N$. Therefore, we have $N\in \mathbb{F}_{\sigma}$ if and only if $N$ does not have elements of order $p$.\\
Thus, by Remark \ref{8}(2) and Example \ref{ExampleZclass}(1), $N\in \mathbb{P}_{\sigma}$ if and only if $N$ has no elements of order p or $N$ is a $p$-group.
\end{itemize}

\end{example}

\begin{example} $\{\alpha_{\mathbb{Z}_p}^{\mathbb{Z}_p} \mid \text{ \textit{p} is prime}\}$-first $\mathbb{Z}$-modules.\\
Let  $\mathscr{A}=\{\alpha_{\mathbb{Z}_p}^{\mathbb{Z}_p} \mid p$ is prime$\}$. The following statements hold:
\begin{enumerate}
    \item 
    There are no non zero  $\mathscr{A}$-fully first $\mathbb{Z}$-modules.
    \item 
    For a non-zero $\mathbb{Z}$-module $M$ we have the following statements that are equivalent:
    \begin{itemize}
        \item[(a)] $M$ is $\mathscr{A}$-first.
        \item[(b)] $M$ has no finite-order elements, or $M$ is a $p$-group for some prime $p$.
    \end{itemize}
\end{enumerate}
\begin{explanation}
\begin{itemize}
    \item[ ] 
\end{itemize}
\begin{enumerate}
\item Let $M \in \mathscr{P}_{\mathscr{A}}$ and $p,q\in \mathbb{Z}$ two different primes. If $M \neq 0$, then $\alpha_{\mathbb{Z}_p}^{\mathbb{Z}_p}(M) \neq 0$. Consequently,  
$\alpha_{\mathbb{Z}_q}^{\mathbb{Z}_q}(\alpha_{\mathbb{Z}_p}^{\mathbb{Z}_p}(M)) \neq 0$,  
which is a contradiction since  
$\alpha_{\mathbb{Z}_q}^{\mathbb{Z}_q} \alpha_{\mathbb{Z}_p}^{\mathbb{Z}_p} = \underline{0}$.  
Therefore, $M = 0$.\\

Note that  $N\in \mathbb{F}_{\mathscr{A}}$ if and only if $N$ has no nontrivial elements of finite order, i.e., $N$ is a torsion-free abelian group.\\

\item Let $M \in \zmod$. 

\begin{itemize}
    \item[$(a) \Rightarrow (b)$]  If $M \in \mathbb{P}_\mathscr{A}$ and $M \notin \mathbb{F}_\mathscr{A}$, then there exists a prime $p \in \mathbb{Z}$ such that $M \notin \mathbb{F}_\sigma$, where $\sigma = \alpha_{\mathbb{Z}_p}^{\mathbb{Z}_p}$. Consequently, we have $\alpha_{\mathbb{Z}_p}^{\mathbb{Z}_p}(M) \neq 0$. This implies that for each $0 \neq N \leq M$, $\alpha_{\mathbb{Z}_p}^{\mathbb{Z}_p}(N) \neq 0$. Hence, each nonzero submodule of $M$ contains an element of order $p$, indicating that $M$ is a $p$-group.\\
 If $M\in \mathbb{F}_\mathscr{A}$ then $M\in \mathbb{P}_\mathscr{A}$ by Remark \ref{8}(4).
 \item[$(b)\Rightarrow(a)$] If $p\in \mathbb{Z}$ is a prime and $M$ is a nonzero $p$-group,  then,   by Lagrange's Theorem, $M$ has no elements of order $q$, for each prime $q\neq p$, so that $\alpha_{\mathbb{Z}_q}^{\mathbb{Z}_q}(M)= 0$.
 If $N$ is a nonzero submodule of $M$, then there exists $x\in N\setminus \{0\}$ with order power of $p$, so, by Cauchy's Theorem, there exists $y\in \mathbb{Z}x$ with order $p$, thus 
 \begin{center}
           \begin{tikzcd}[row sep=tiny]
        f:\mathbb{Z}_p\arrow[r]& N \\
                a+p\mathbb{Z}\arrow[r, mapsto] & ay,
           \end{tikzcd}
        \end{center}
 $f$ is a well defined nonzero $\mathbb{Z}$-morphism, which implies $\alpha_{\mathbb{Z}_p}^{\mathbb{Z}_p}(N)\neq 0$. Then, $M$ is an $\mathscr{A}$-first module. 
\end{itemize}

\end{enumerate}
\end{explanation}
\end{example}

\begin{proposition}\label{12}
For a ring $R$, we have $\mathscr{P}_{soc}=\rmod$ if and only if  $\mathbb{P}_{soc}=\rmod$.
\end{proposition}

\begin{proof}
By Remark \ref{8},  $2$, it is clear that if $\mathscr{P}_{soc}=\rmod$,  then $\mathbb{P}_{soc}=\rmod$. Now, if $\mathbb{P}_{soc}=\rmod$ and $N\in R$-Mod$\setminus \{ 0 \}$, $N$ embeds in $E(C_0)^X$ where $C_0=\bigoplus\limits_{S\in \rsimp}S$,  for a set $X$. Also, $soc(E(C_0)^X)\neq 0$ and $E(C_0)^X\in \mathbb{P}_{soc} $, then $E(C_0)^X\in \mathscr{P}_{soc} $, so $N\in \mathscr{P}_{soc} $. Therefore $\mathscr{P}_{soc}=\rmod$.
\end{proof}

\begin{remark}\label{13}
If $N\in R$-Mod is   semisimple, then $N\in \mathbb{P}_{\rlep}$ if and only if $N$ is  semisimple homogeneous. This is because as  $\alpha_S^S$ is a left exact preradical for any $S\in \rsimp$, we can adapt the proof of Lemma  \ref{5} to $\rlep$.
\end{remark}

\begin{proposition}\label{8.5}
$\mathscr{P}_{soc}=\rmod$ if and only if $R$ is a left semiartinian ring.
\end{proposition}

\begin{proof}
 $R$ is semiartinian if and only if $soc(M)\neq 0$ for all $0\neq M\in \rmod$.  From this it  follows that $\mathscr{P}_{soc}=\rmod$ if and only if $R$ is left semiartinian.
\end{proof}

\begin{remark}
If $\sigma,\tau\in \rlep$, then $\sigma(\tau(M))=\sigma(M)\cap \tau(M)=\tau(\sigma(M))$ for all $M\in \rmod$. So $\sigma \circ\tau=\tau\circ\sigma$ if $\sigma$ and $\tau$ are left exact preradicals.
\end{remark}

\begin{theorem}\label{14}
For a ring $R$,  the following statements are equivalent:
\begin{itemize}
    \item[(1)] $R$ is a left semiartinian left local ring,
    \item[(2)] $ \mathbb{P}_{\rlep}=\rmod$.
\end{itemize}

\end{theorem}

\begin{proof}
\begin{itemize}
    \item[]
     \item[(1)$\Rightarrow$ (2)] Let $M\in \rmod$, $0\neq N\leq M$ and let  $\sigma\in \rlep$ such that $\sigma(M)\neq 0$. As $R$ is left semiartinian, then $soc(\sigma(M))\neq 0$ and $soc(N)\neq 0$. Besides, as $\sigma$  is left exact, then $\sigma\circ soc=soc\circ\sigma$, so $\sigma(soc(M))\neq 0$.\\
     On the other hand, as  $R$ is left local, then  $\sigma(soc(M))\neq 0$  implies that $\sigma(S)=S$  embeds in  $N$, where $S$ is a simple module. Then $\sigma(S)$ embeds in $\sigma(N)$, so $\sigma(N)\neq 0$. Therefore, $M\in \mathbb{P}_{\rlep}$.
    \item[(2)$\Rightarrow$ (1)] By Remark \ref{8}, $1$, we have that, $ \mathbb{P}_{soc}=\rmod$. So  from Remark \ref{8.5} and Proposition \ref{12}, it follows that $R$ is a left semiartinian ring. On the other hand, from Remark \ref{13}, it follows that $R$ is left local.
\end{itemize}
\end{proof}

\section{Rings for which every module is $\rpr$-first.}

\begin{proposition}\label{14.1}
Let $S\in \rsimp$ such that $S\lneq E(S)$. Then $S$ is superfluous in  $E(S)$.
\end{proposition}

\begin{proof}
Let $N\leq E(S)$ such that  $S+N=E(S)$. Since $S\lneq E(S)$, we have $N\neq 0$. Moreover, $S$  is essential  $E(S)$, so $0\neq S\cap N\leq S $.  On the other hand, $S$ is simple, then $S=S\cap N$, which implies that  $S\leq N$ and $N=S+N=E(S)$. 
\end{proof}

\begin{remark}\label{14.2}
 $ \mathbb{P}_{\rpr}=\rmod$ implies $ \mathbb{P}_{\rlep}=\rmod$. This because, by Remark \ref{8}, (5), $\mathbb{P}_{\rpr}\subseteq\mathbb{P}_{\rlep}$. 
\end{remark}

\begin{theorem}\label{14.3}
For a ring $R$,  the following statements are equivalent:
\begin{itemize}
    \item[(1)] $R$ is  a left semiartinian left local V-ring.
    \item[(2)] $ \mathbb{P}_{\rpr}=\rmod$.
     \item[(3)] $R$ is a semisimple homogeneous ring.
\end{itemize}

\end{theorem}

\begin{proof}
\begin{itemize}
    \item[]
     \item[(1)$\Rightarrow$ (2)] Let $\alpha\in \rpr\setminus \{0\}$. Then there exists $S\in \rsimp$ such that  $\alpha_S^{E(S)}\preceq \alpha$,  but $R$ is   a left local $V$-ring, which implies that  $soc=\alpha_S^S=\alpha_S^{E(S)}$. As $R$ is semiartinian, then for all  nonzero $M\in R$-Mod,  we have that: \begin{center}
         $0\neq soc(M)=\alpha_S^S(M)=\alpha_S^{E(S)}(M)\leq \alpha (M)
     $.
     \end{center} It follows that  $\mathbb{P}_{\rpr}=\rmod$.
    \item[(2)$\Rightarrow$ (1)]By Remark  \ref{14.2},  we have that  $\mathbb{P}_{\rlep}=\rmod$, which implies, by Theorem \ref{14},  that  $R$ is left  semiartinian left local.\\
    Let $S\in \rsimp$.  If  $S\lneq E(S)$, then $S$ is superfluous in $E(S)$ by Proposition \ref{14.1},  which implies that  $rad(E(S))\neq 0$ y $rad(S)=0$, so $E(S)\notin \mathbb{P}_{\rpr}$,  which is a contradiction.\\
   Therefore,  $S=E(S)$  for all $S\in \rsimp$, which implies that $R$ is a V-ring.
   \item[(3)$\Rightarrow$ (1)] We already know that $R$ is left local, and that every R-module is injective and projective since $R$ is semisimple. In particular, every $R$-simple module is injective, so $R$ is a $V$-ring.\\
   It is clear that $R$ is left semiartinian.
   \item[(1)$\Rightarrow$ (3)] Since $R$ is left semiartinian, there exists a simple submodule  $S$  of $R$. Furthermore, since $R$ is a $V$-ring, then $S$ is a direct summand of $R$, so it is projective. 
Suppose $soc (R)\lneq R$. Then there is a maximal ideal $M$ of $R$ such that  $soc(R)\leq M$. Moreover, since $R$ is left local,  we have that $R/M\cong S$,  so the following sequence splits:\\
\begin{tikzcd}
 & & & 0\arrow{r}& M \arrow{r}& R \ar{r} & R/M\ar{r}& 0.
\end{tikzcd} \\
Then there exists a submodule $K$ of $R$ such that $K\cap M=\{ 0 \}$, 
and $K\cong R/M$, so $K\cap soc(R)=\{ 0\}$ and $K$ is simple,  which is a contradiction. Therefore, $R=soc(R)$ implies $R$ es semisimple and, by hypothesis, left local.
\end{itemize}
\end{proof}

\begin{definition}
Let $R$ be a  ring. $R$ is a BKN ring if  $Hom_R(M,N)\neq 0$ for any nonzero $M,N\in \rmod$, see Proposition VI.2.3 of \cite{Bic}.
\end{definition}

\begin{proposition}\label{error1}
Let $R$ be a ring. If $\mathbb{P}_{\rpr}=\rmod$, then $R$ is a BKN ring.
\end{proposition}

\begin{proof}
Theorem \ref{14.3} shows that  $R$ is a left local semiartinian $V$-ring. Let  $M, N\in \rmod$ be nonzero, and assume they are the only $R$-simple module except for  isomorphic copies. Since $R$ is semiartinian,  there exists $0\neq f\in Hom_R(S,N)$. Moreover, since $S$ is simple, we have that $f$ is a monomorphism. On the other hand, since $R$ is left local and $V$-ring, we have  $0\neq g\in Hom_R(M,S)$. Then $f\circ g\in Hom_R(M,N)$ and $f\circ g\neq 0$ since $0\neq f$ is a monomorphism, and $g\neq 0$. Therefore, $R$ is a BKN ring.
\end{proof}

\begin{example}
Consider the ring  $\mathbb{Z}_{p^2}$ with $p$ a prime.   $\mathbb{Z}_{p^2}$ is a left local semiartinian MAX ring,  a BKN ring. Moreover, we have that $E(\mathbb{Z}_p)=\mathbb{Z}_{p^2}$, so  $\mathbb{Z}_p$ is not a  $V$-ring. Then, by Theorem  \ref{14.3}, we have  $\mathbb{P}_{\rpr}\neq \rmod$, which shows that Proposition \ref{error1} does not hold.
\end{example}

\noindent\textbf{Luis Fernando Garc\'ia-Mora}\\
Departamento de Matemáticas, Facultad de Ciencias,  \\ 
Universidad Nacional Autónoma de México, Mexico City, Mexico. \\ 
\textbf{e-mail:} \textit{lu1sgarc1agm1995@gmail.com}

\noindent\textbf{Hugo Alberto Rinc\'on-Mej\'ia}\\
Departamento de Matemáticas, Facultad de Ciencias,  \\ 
Universidad Nacional Autónoma de México, Mexico City, Mexico. \\ 
\textbf{e-mail:} \textit{hurincon@gmail.com}\\

\end{document}